\documentclass{article}
\usepackage{graphicx} 

\usepackage{graphicx}          
\usepackage{amsmath,amssymb}   
\usepackage{amsthm}            
\usepackage{booktabs}          
\usepackage{enumitem}          
\usepackage{float}
\usepackage{tikz}

\usepackage{cite}
\usepackage{caption}          
\usepackage{authblk}   

\usetikzlibrary{positioning}
\usetikzlibrary{fit}
\usetikzlibrary{shapes.geometric}

\newtheoremstyle{myplain}%
  {3pt}
  {3pt}
  {\itshape}
  {}
  {\bfseries}
  {}
  { }
  {\thmname{#1}\thmnumber{ #2}\thmnote{ (#3)}}

\theoremstyle{myplain}
\newtheorem{theorem}{Theorem}[section]
\newtheorem{lemma}[theorem]{Lemma}
\newtheorem{proposition}[theorem]{Proposition}
\newtheorem{corollary}[theorem]{Corollary}
\newtheorem{example}[theorem]{Example}

\newtheorem{definition}[theorem]{Definition}

\title{\textbf{Divisibility graphs of semigroups as weighted poset blow-ups
}}
\author[1]{\textbf{Xuan Yang}}
\author[1]{\textbf{Yong Shao$^{\dagger}$}}   
\author[1]{\textbf{Lingli Zeng}}

\affil[1]{School of Mathematics, Northwest University, Xi'an 710127, Shaanxi, People's Republic of China}

\date{}
\begin{document}
\maketitle


\renewcommand{\thefootnote}{\textdagger}
\footnotetext[1]{Corresponding author.\\ \hspace*{1em} E-mail addresses: xuanyangmath@163.com (X.Yang), yongshaomath@126.com (Y.Shao), zengll929@nwu.edu.cn (L.L.Zeng).}

\begin{abstract}
We develop a structural theory of divisibility graphs of arbitrary semigroups. 
For a semigroup $S$, let \(S/\mathcal J\) be the poset of Green’s \(\mathcal J\)-classes ordered by inclusion of principal ideals, and equip each \(\mathcal J\)-class \(J\) with the weight \(|J|\). 
We prove that the divisibility graph \(\operatorname{Div}(S)\) is obtained from the comparability directed graph $Comp(S/\mathcal{J})$ by replacing each \(\mathcal{J}\)-class with a complete symmetric graph of order \(|J|\) and by replacing each directed edge $J_u\to J_v$ in $Comp(S/\mathcal{J})$ with all directed edges from every element of $J_u$ to every element of $J_v$.
Consequently, \(\operatorname{Div}(S)\) is determined, up to isomorphism, by the weighted poset \((S/\mathcal{J},\{J\}_{J\in S/\mathcal{J}})\). This yields structural descriptions of its strongly connected components, maximal unilaterally connected subgraphs.
The main contribution is the reduction of divisibility graph problems to weighted quotient-poset problems. As applications, we obtain reconstruction and isomorphism criteria and prove an inheritance theorem for divisibility \(D\)-saturated subsemigroups.

\end{abstract}

\textbf{Keywords:} Divisibility graph; Semigroup; Green’s \(\mathcal J\)-relation; Quotient poset; \(D\)-saturated property.

\section{Introduction}
Graph structures associated with semigroups have been widely studied in recent years, including Cayley graphs, power graphs, commuting graphs (see~\cite{Kelarev2006, Luo2011, Khosravi2016, AsheghBonabi2022, Cheng2024, Cheng2025, Zeng2025, Araujo2015, Paulista2025, KelarevPraeger2003, Khosravi2014, Wang2017, Yang2026, Dalal2024, Araujo2011}), and divisibility graphs.
Among them, divisibility graphs are closely related to the ideal structure of a semigroup and to Green’s relations. 
For a semigroup $S$, the divisibility graph $\operatorname{Div}(S)$ is a directed graph whose vertex set is $S$, and for all $u,v\in S$ with $u\ne v$, $u\to v$ if and only if $u\in S^1vS^1$ (see~\cite{Kelarev2002}). 
Thus, $\operatorname{Div}(S)$ naturally reflects the principal ideal structure of $S$ and provides a bridge between semigroup theory and algebraic graph theory.

Kelarev and Quinn~\cite{Kelarev2000} first introduced the $D$-saturated property in the context of power graphs of groups.  
Subsequently, the investigation of the $D$-saturated property was extended to semigroups and to several other graph classes, including power graphs, annihilator graphs, Cayley graphs, and divisibility graphs (see~\cite{Kelarev2002, Kelarev2003, Yang2010, Khosravi2015}).
In particular, for divisibility graphs, a series of studies has been devoted to characterizing the $D$-saturated property for various classes of semigroups, including commutative semigroups, monomial matrix semigroups, and 
completely regular semigroups (see~\cite{Kelarev2002, Quinn2004, Hao2010}).  
More recently, Guo et al.~\cite{Guo2025} established a key reduction theorem for arbitrary semigroups, stating that a semigroup is divisibility $D$-saturated if and only if the quotient poset $(S/\mathcal{J},\le)$ is $D$-saturated.
This result shows that $(S/\mathcal{J},\le)$ is an important tool for studying divisibility graphs.
It also raises a stronger question: can the entire graph $\operatorname{Div}(S)$ be recovered from $(S/\mathcal{J},\le)$ together with the sizes of its $\mathcal{J}$-classes?

In this paper, rather than addressing only the $D$-saturated property, we give a complete structural reconstruction of $\operatorname{Div}(S)$ by showing that, for an arbitrary semigroup $S$, its divisibility graph is determined up to isomorphism by the weighted quotient poset $(S/\mathcal{J}, \{J\}_{J\in S/\mathcal{J}})$.
This description immediately yields several consequences.
In particular, the strongly connected components of $\operatorname{Div}(S)$ are precisely the $\mathcal{J}$-classes, the maximal unilaterally connected subgraphs correspond to the maximal chains in $(S/\mathcal{J},\le)$, and $\operatorname{Div}(S)$ is always weakly connected, with its underlying undirected graph having diameter at most $2$.
Furthermore, we obtain an isomorphism criterion for the divisibility graphs of semigroups.

As further applications, we introduce the notion of divisibility-preserving subsemigroups and prove that this condition ensures the inheritance of the divisibility $D$-saturated property from a semigroup to its subsemigroups.
The paper is organized as follows.  
Section~$2$ recalls the necessary preliminaries.  
Section~$3$ establishes the blow-up theorem and its structural consequences.  
Section~$4$ studies divisibility $D$-saturated subsemigroups and proves the inheritance theorem.

\section{Preliminaries}
For a semigroup $S$ and $a\in S$, we have that $S^1aS^1$ is the principal ideal of $S$ generated by $a$, where $S^1$ is $S$ adjoining an identity if necessary. An equivalence relation $\mathcal{J}$ on $S$ is defined by the rule that $a\mathcal{J}b$ if and only if $S^1aS^1=S^1bS^1$. 
The $\mathcal{J}$-class containing the element $a$ is denoted by $J_a$. 
The quotient set $S/{\mathcal{J}}$ is the set of all $\mathcal{J}$-classes, i.e., $S/{\mathcal{J}} = \{ J_a \mid a\in S \}$. On the quotient set $S/{\mathcal{J}}$ we define a partial order $\le$ by $J_a \le J_b \iff S^1aS^1\subseteq S^1bS^1$ (see~\cite{Howie1995}, Section 2.1).
Naturally, $(S/{\mathcal{J}},\le)$ is a poset.

For an undirected graph $G$, we denote its vertex set by $V(G)$ and its edge set by $E(G)$. 
A path in $G$ is a sequence of distinct vertices $v_0, v_1, \dots, v_k$ ($k>0$) such that $\{v_i, v_{i+1}\}\in E(G)$ for each $i=0,\dots,k-1$.  
For two distinct vertices $u,v\in V(G)$, we say that $u$ and $v$ are connected if 
there is a path in $G$ joining $u$ and $v$.
The graph $G$ is called connected if every pair of distinct vertices is connected.
A connected component of $G$ is a maximal connected induced subgraph of $G$ 
(see~\cite{Bondy2008}).

Similarly, for a directed graph $D$, we denote its vertex set by $V(D)$ and its directed edge set by $E(D)$. 
Throughout this paper, we shall write $(u,v)\in E(D)$ as $u\rightarrow v$.
A directed path in $D$ is a sequence of distinct vertices $v_0, v_1, \dots, v_k (k>0)$ such that $v_i\to v_{i+1}$ for each $i = 0,\dots,k-1$. 
A directed cycle in $D$ is a sequence of vertices $v_0, v_1, \dots, v_k(k>0)$ such that $v_0 \to v_1 \to \cdots \to v_k$ with $v_0 = v_k$, and all vertices $v_0, v_1, \dots, v_{k-1}$ are distinct. 
For two distinct vertices $u,v\in V(D)$ we say that $u$ and $v$ are strongly connected if there is a directed path from $u$ to $v$ and also one from $v$ to $u$; $u$ and $v$ are unilaterally connected if there is a directed path in at least one direction; and they are weakly connected if they are connected in the underlying undirected graph obtained by ignoring edge orientations. The graph $D$ is called strongly (respectively, weakly, unilaterally) connected if every pair of distinct vertices is strongly (respectively, weakly, unilaterally) connected; a strongly (respectively, weakly) connected component is a maximal induced subgraph of $D$ that is strongly (respectively, weakly) connected (see~\cite{Bondy2008}).
$D$ is called acyclic if it contains no directed cycles.
The underlying undirected graph of a directed graph $D$ is the undirected 
graph with vertex set $V(D)$ in which two distinct vertices $u$ and $v$ are 
adjacent if and only if $u\to v$ or $v\to u$ is a directed edge of $D$.  
Obviously, two vertices are weakly connected in $D$ precisely when they are connected in the underlying undirected graph of $D$.

Let $S$ be a semigroup and $T$ be a subset of $S$, we shall say that $T$ is strongly (respectively, weakly, unilaterally) connected if the subgraph of $\operatorname{Div}(S)$ induced by $T$, denoted by $\Gamma(T)$, is strongly (respectively, weakly, unilaterally) connected.

Let $(P,\le)$ be a poset. The comparability directed graph of $P$, denoted by $Comp(P)$, is a directed graph whose vertex set is $P$, and for all $x,y\in P$ with $x\ne y$, $x\to y$ if and only if $x<y$ (see~\cite{Bier2018}).
Let $A$ be a subset of $P$. $A$ is called an anti-chain if every two distinct elements of $A$ are incomparable.

Let $D$ be a finite graph. A semigroup $S$ is called divisibility $D$-saturated if for every infinite subset $T$ of $S$, the divisibility graph $\operatorname{Div}(S)$ contains a subgraph isomorphic to $D$ whose vertex set is entirely contained in $T$ (see~\cite{Kelarev2002}). Clearly, any finite semigroup is divisibility $D$-saturated. Recall that a null graph is a graph without edges. It is trivial that any semigroup is divisibility $D$-saturated if $D$ is a null graph.

In this paper, $\mathbb{N}$ denotes the set of all positive integers.
A complete symmetric graph on a vertex set $V$ is a directed graph in which, for every pair of distinct vertices $u,v\in V$, both directed edges $u\to v$ and $v\to u$ exist.

We refer the reader to~\cite{Howie1995} and~\cite{Bondy2008} for undefined terms and notation, and we review some known results that will be used frequently in the following.

\begin{lemma}
\textup{\cite{Howie1995}}
\label{lemma: simple-equivalent}
Let $S$ be a semigroup. Then $S$ is simple if and only if $\mathcal{J}=S\times S$. 
\end{lemma}

\begin{lemma}[Zorn's Lemma]
\textup{\cite{Dummit1999}}
\label{lemma: Zorn}
Let $(P,\le)$ be a nonempty poset. If every chain in $P$ has an upper bound in $P$, then $P$ contains a maximal element.   
\end{lemma}

By Zorn's Lemma, every element of a poset is contained in a 
maximal chain. Indeed, assume that $(P,\le)$ is a nonempty poset and $x\in P$, let $\mathcal{F}$ be the family of all chains in $P$ that contain $x$, ordered by inclusion. Then $\mathcal{F}$ is nonempty since $\{x\}\in\mathcal{F}$, and $\mathcal{F}$ satisfies the hypotheses of Zorn's Lemma. Thus $\mathcal{F}$ has a maximal element $M$, which is a maximal chain in $P$ containing $x$.

\section{The structure of the divisibility graphs}
In this section, we explore the fundamental properties and connectivity of divisibility graphs of semigroups and demonstrate the structure of divisibility graphs.

By the definition of a divisibility graph, if $u\rightarrow v$ in $\operatorname{Div}(S)$, then $u \in S^1 v S^1$ and so $S^1 u S^1 \subseteq S^1 v S^1$. Thus we obtain the following lemma.

\begin{lemma}
\label{lemma: edge of Div(S)-equivalent}
Let $S$ be a semigroup and $u,v\in S$ with $u\ne v$. Then the following statements are equivalent.
\begin{enumerate}[label=(\roman*)]
  \item $u\rightarrow v$ in $\operatorname{Div}(S)$;
  \item $S^1 u S^1 \subseteq S^1 v S^1$;
  \item $J_u \le J_v$.
\end{enumerate}
\end{lemma}

By Lemma \ref{lemma: edge of Div(S)-equivalent}, we have the following proposition and corollary.

\begin{proposition}
\label{prop: Divisibility graph structure}
Let $S$ be a semigroup and $u,v\in S$ with $u\ne v$. Then
\begin{enumerate}[label=(\roman*)]
\item If $J_u=J_v$, then both $u\to v$ and $v\to u$ occur in $\operatorname{Div}(S)$;
\item If $J_u<J_v$, then every element of $J_u$ has a directed edge to every element of $J_v$, and no edge points in the opposite direction;
\item If $J_u$ and $J_v$ are incomparable, then no edge occurs between these two classes. 
\end{enumerate} 
\end{proposition}

\begin{corollary}
\label{coro: Div(S) two edge----J, transitivity, directed cycle in J}
Let $S$ be a semigroup and $u,v,w\in S$. Then
\begin{enumerate}[label=(\roman*)]
 \item For $u\ne v$, both $u\rightarrow v$ and $v\rightarrow u$ occur in $\operatorname{Div}(S)$ if and only if $u\mathcal{J}v$;
  \item If $u\rightarrow v$ and $v\rightarrow w$ in $\operatorname{Div}(S)$, then $u\rightarrow w$ in $\operatorname{Div}(S)$ whenever $u\ne w$;
  \item Any directed cycle in $\operatorname{Div}(S)$ must be completely contained within some $\mathcal{J}$-class of $S$.
\end{enumerate} 
\end{corollary}

In Definition~\ref{def: Weighted blow-up} and Theorem~\ref{thm: Divisibility graph blow-up theorem}, a poset $(P,\le)$ is simply denoted by $P$.

\begin{definition}[Weighted blow-up]
\label{def: Weighted blow-up}
Let $P$ be a poset. For each $x\in P$, let $I_x$ be a nonempty set. 
The order $|I_x|$ is called the weight of $x$. 
We call the pair $(P,\{I_x\}_{x\in P})$ a weighted poset.

The weighted blow-up of $(P,\{I_x\}_{x\in P})$, denoted by 
$\mathsf{Bl}(P,\{I_x\}_{x\in P})$, is a directed graph whose vertex set is 
$V=\{(x,\alpha) \mid x\in P,\ \alpha\in I_x\}$,
and for all $(x,\alpha),(y,\beta)\in V$ with $(x,\alpha)\ne (y,\beta)$,
$(x,\alpha)\to(y,\beta)$ if and only if $x\le y$.
\end{definition}

Viewed in this way, the weighted blow-up construction is the core model of this section: each $\mathcal{J}$-class plays the role of a fiber, while the quotient poset records the directed relations between different fibers. 
The following theorem makes this correspondence precise.

\begin{theorem}[Divisibility graph blow-up theorem]
\label{thm: Divisibility graph blow-up theorem}
Let $S$ be a semigroup. Then
$\operatorname{Div}(S) \cong \mathsf{Bl}(S/\mathcal{J},\{J\}_{J\in S/\mathcal{J}})$.
\end{theorem}

In the following, we investigate the strongly connected components, maximal unilaterally connected subgraphs, and weak connectivity of $\operatorname{Div}(S)$.

In Proposition \ref{prop:connectivity(strongly,unila,weak)}, we identify a subset $U\subseteq S$ with the subgraph of $\operatorname{Div}(S)$ induced by $U$.

\begin{proposition}
\label{prop:connectivity(strongly,unila,weak)}
Let $S$ be a semigroup. Then
\begin{enumerate}[label=(\roman*)]
    \item The strongly connected components of $\operatorname{Div}(S)$ are precisely the $\mathcal{J}$-classes of $S$.
    \item The maximal unilaterally connected subgraphs of $\operatorname{Div}(S)$ are precisely the sets $\bigcup_{J\in M} J$, where $M$ is a maximal chain in
          $(S/{\mathcal{J}},\le)$.
    \item $\operatorname{Div}(S)$ is always weakly connected. Furthermore, the underlying undirected graph of $\operatorname{Div}(S)$ has diameter 
at most $2$.    
\end{enumerate}
\end{proposition}

\begin{proof}
(i) If $U$ is a strongly connected component of $\operatorname{Div}(S)$,
then for any distinct $u,v\in U$, there exist directed paths from $u$ to $v$ and from $v$ to $u$, i.e., $u \rightarrow c_1 \rightarrow c_2 \rightarrow \cdots \rightarrow c_k \rightarrow v $ and $v \rightarrow d_1 \rightarrow d_2 \rightarrow \cdots \rightarrow d_l \rightarrow u$. 
By Corollary \ref{coro: Div(S) two edge----J, transitivity, directed cycle in J}(ii), we have $u \rightarrow v$ and $v \rightarrow u$. 
Corollary \ref{coro: Div(S) two edge----J, transitivity, directed cycle in J}(i) implies that $u\mathcal{J}v$.
Hence $U\subseteq J_u$.
Obviously, $J_u$ is strongly connected and so $U=J_u$. 
Conversely, let $J$ be a $\mathcal{J}$-class of $S$, which is strongly connected.
By Proposition~\ref{prop: Divisibility graph structure}(ii)(iii), between different $\mathcal{J}$-classes, all edges are unidirectional.
Hence $J$ is a strongly connected component.

(ii) If a subset $U\subseteq S$ is maximal unilaterally connected,
then for any distinct  $u,v\in U$, there exists a directed path from $u$ to $v$ or from $v$ to $u$. 
By Corollary \ref{coro: Div(S) two edge----J, transitivity, directed cycle in J}(ii), we have $u \rightarrow v$ or $v \rightarrow u$
and so $J_u$ and $J_v$ are comparable.
Hence the set $M_U=\{ J_u \mid u \in U\}$ is a chain in $(S/\mathcal{J}, \le)$.
Suppose that there exists $w \in S \setminus U$ such that $M_U\cup \{J_w\}$ is a chain in ($S/\mathcal{J}$, $\le$). 
Then for any $u\in U$, $J_u$ and $J_w$ are comparable, and so $u \rightarrow w$ or $w \rightarrow u$. 
Hence $w \in U$, contradicting the assumption that $w \in S \setminus U$. 
Thus, $M_U$ is a maximal chain in ($S/\mathcal{J}$, $\le$) and so $U\subseteq\bigcup_{J\in M_U} J$.
Conversely, if $x\in J$ for some $J\in M_U$, then $J=J_u$ for some $u\in U$, so
$x\in J_u$. 
Hence $x\in U$ since $U$ is maximal unilaterally connected.  Thus $\bigcup_{J\in M_U} J \subseteq U$. Therefore $U=\bigcup_{J\in M_U} J$.

Conversely, let $M$ be a maximal chain in $(S/{\mathcal{J}},\le)$ and let
$U_M = \bigcup_{J\in M} J$.  
For any distinct $u,v\in U_M$, $J_u$ and $J_v$ are comparable, so $u\to v$ or $v\to u$. 
Hence $U_M$ is unilaterally connected.
Suppose that there exists $w\in S\setminus U_M$ such that
$U_M\cup\{w\}$ is unilaterally connected.  
Then $J_w\notin M$ and for any $u\in U_M$, there exists a directed path from $u$ to $w$ or from $w$ to $u$, so $u\to w$ or $w\to u$, which implies that $J_u$ and $J_w$ are comparable.
Therefore $M\cup\{J_w\}$ is a chain, contradicting the maximality of $M$.  
Hence $U_M$ is maximal.

(iii) For any distinct $u,v\in S$, let $w = uv$.  Then $w=uv\in S^1uS^1\cap S^1vS^1$.
Hence $J_w \le J_u$ and $J_w \le J_v$ and so $w\to u$ and $w\to v$ whenever 
$w\ne u$ and $w\ne v$, respectively.

If $w\ne u$ and $w\ne v$, then $u\leftarrow  w\to v$, so the distance between $u$ and $v$ in the underlying undirected graph is at most $2$.

If $w = u$, then $u = uv \in S^1 v S^1$.  Since $u\ne v$, this gives the 
directed edge $u\to v$, so their distance in the underlying undirected graph is exactly $1$.

Similarly, if $w = v$, then $v\to u$ is a directed edge, so their distance in the 
underlying undirected graph is exactly $1$.

Therefore, any two distinct vertices of $\operatorname{Div}(S)$ are connected in 
the underlying undirected graph by a path of length at most $2$. This shows that 
the underlying undirected graph is connected, and that its diameter is at most 
$2$. Consequently, $\operatorname{Div}(S)$ is weakly connected.

\end{proof}

By Proposition \ref{prop:connectivity(strongly,unila,weak)} and Corollary \ref{coro: Div(S) two edge----J, transitivity, directed cycle in J}, we have that

\begin{corollary}
\label{coro: Div(S) is strongly connected----equivqlent }
Let $S$ be a semigroup. Then the following statements are equivalent.
\begin{enumerate}[label=(\roman*)]
  \item $\operatorname{Div}(S)$ is a complete symmetric graph;
  \item $\operatorname{Div}(S)$ is strongly connected;
  \item $\mathcal{J}=\nabla$;
  \item S is simple.
\end{enumerate}
\end{corollary}

A poset is totally ordered if and only if it has exactly one maximal chain. 
Indeed, if $P$ is totally ordered, then $P$ itself is a chain. 
If $M$ is any maximal chain, then $M\subseteq P$. 
If $M\ne P$, choose $x\in P\setminus M$. 
Since $P$ is totally ordered, $x$ is comparable with every element of $M$, so 
$M\cup\{x\}$ is still a chain, which contradicts the maximality of $M$. 
Hence $M=P$, and $P$ is the unique maximal chain.
Conversely, assume that $P$ has exactly one maximal chain $M$. 
For any $x\in P$, $\{x\}$ is a chain. 
By Lemma \ref{lemma: Zorn}, $\{x\}$ can be extended to a maximal chain, which 
must be $M$. 
Thus $x\in M$. 
Hence $P=M$, so $P$ is a chain and therefore totally ordered.

Thus we immediately obtain the following corollary.

\begin{corollary}
\label{coro: Div(S) is unilaterally connected-----equivalent }
Let $S$ be a semigroup. Then $\operatorname{Div}(S)$ is unilaterally connected if and only if ($S/\mathcal{J}$, $\le$) is a totally ordered set.
\end{corollary}

\begin{proposition}
\label{prop: maximal chains}
Let $S$ be a semigroup and $\mathcal{M}$ be the set of all maximal chains in $(S/\mathcal{J},\le)$. Then for all $M,N\in\mathcal{M}$, there exists a finite sequence $M=M_0,M_1,\cdots,M_k=N$ in $\mathcal{M}$ such that for each $i=0,1,\cdots,k-1$, $M_i\cap M_{i+1}\ne\emptyset$.  
\end{proposition}

\begin{proof}
Define a graph $D(\mathcal{M})$ on $\mathcal{M}$ whose vertex set is $\mathcal{M}$, where two distinct maximal chains $M_1$ and $M_2$ are adjacent if and only if $M_1\cap M_2\ne\emptyset$. 

We claim that $D(\mathcal{M})$ is a connected graph. If $D(\mathcal{M})$ is not connected, and let $\mathcal{N}_1,\mathcal{N}_2,\dots$ 
be its connected components.  Take $\mathcal{M}_1=\mathcal{N}_1$ and let $\mathcal{M}_2$ be the union of all the remaining components.
Then $\mathcal{M}$ can be divided into two nonempty, disjoint subsets $\mathcal{M}_1$ and $\mathcal{M}_2$ such that for any $M\in\mathcal{M}_1$ and $N\in\mathcal{M}_2$, $M\cap N=\emptyset$. 
Let $C_i=\bigcup_{M\in\mathcal{M}_i} M, i=1,2$. Then $S/\mathcal{J}=C_1\cup C_2$ and $C_1\cap C_2=\emptyset$. We claim that there is no directed edge between $C_1$ and $C_2$ in $Comp(S/\mathcal{J})$. If not, there exist $J_x\in C_1$ and $J_y\in C_2$ such that $J_x\to J_y$ or $J_y\to J_x$, which means that $J_x$ and $J_y$ are comparable. Without loss of generality, assume that $J_x< J_y$. By Lemma~\ref{lemma: Zorn}, the chain $\{J_x,J_y\}$ can be extended to a maximal chain. Therefore, there exists $\overline{M}\in\mathcal{M}$ such that $J_x,J_y\in\overline{M}$. This implies that $C_1\cap\overline{M}\ne\emptyset$ and $C_2\cap\overline{M}\ne\emptyset$. Since $\overline{M}\in\mathcal{M}$, either $\overline{M}\in\mathcal{M}_1$ or $\overline{M}\in\mathcal{M}_2$, which contradicts $C_1\cap\overline{M}\ne\emptyset$ and $C_2\cap\overline{M}\ne\emptyset$. It follows that there is no directed edge between $C_1$ and $C_2$ in $Comp(S/\mathcal{J})$, which contradicts the fact that $Comp(S/\mathcal{J})$ is a weakly connected graph (the proof follows by the same argument as in Proposition~\ref{prop:connectivity(strongly,unila,weak)}(iii)). Therefore, $D(\mathcal{M})$ is a connected graph. This means that for each $M,N\in\mathcal{M}$, there exists a finite sequence $M=M_0,M_1,\cdots,M_k=N$ such that for each $i=0,1,\cdots,k-1$, $M_i\cap M_{i+1}\ne\emptyset$.
\end{proof}

\begin{example}
Let $S=\{0,a,b,c,d,e,f\}$ and the Cayley table is shown in Table \ref{tab:clifford}. 
It is easy to check that $S$ is a Clifford semigroup with identity element $e$. $S0S=\{0\}$, $SaS=SbS=\{a,b,0\}$, $ScS=SdS=\{c,d,0\}$, $SeS=SfS=S$. $S$ has only four $\mathcal{J}$-classes, namely $J_0=\{0\},J_a=\{a,b\},J_c=\{c,d\},J_e=\{e,f\}$. 
\end{example}

The comparability directed graph $Comp(S/\mathcal{J})$ is shown in Figure \ref{fig:diamond(lingxing) Comp(S/J)}. The red maximal chain $M_1=\{J_0,J_a,J_e\}$ and the blue one $M_2=\{J_0,J_c,J_e\}$ cover the entire poset and share $J_0$ and $J_e$.

\begin{figure}[H]   
\centering
\begin{minipage}{0.48\textwidth}
\centering
\captionsetup{skip=0pt}   
\captionof{table}{Cayley table of semigroup $S$}   
\label{tab:clifford}
\renewcommand{\arraystretch}{0.7}   
\setlength{\arraycolsep}{3pt}       
\[
\begin{array}{c|ccccccc}
\cdot & 0 & a & b & c & d & e & f \\ \hline
0     & 0 & 0 & 0 & 0 & 0 & 0 & 0 \\
a     & 0 & a & b & 0 & 0 & a & b \\
b     & 0 & b & a & 0 & 0 & b & a \\
c     & 0 & 0 & 0 & c & d & c & d \\
d     & 0 & 0 & 0 & d & c & d & c \\
e     & 0 & a & b & c & d & e & f \\
f     & 0 & b & a & d & c & f & e
\end{array}
\]
\end{minipage}
\hfill
\begin{minipage}{0.48\textwidth}
\centering
\begin{tikzpicture}[
    scale=0.4,   
    node distance=0.3cm and 1.5cm,
    class/.style={circle, draw, inner sep=1pt, minimum size=0.8cm},
    chain1/.style={red, line width=1.5pt, opacity=0.7},
    chain2/.style={blue, line width=1.5pt, opacity=0.7},
    >=stealth,   
]
\node[class] (J0) {\(J_0\)};
\node[class, above left=of J0] (Ja) {\(J_a\)};
\node[class, above right=of J0] (Jc) {\(J_c\)};
\node[class, above=of J0, yshift=0.6cm] (Je) {\(J_e\)};

\draw[black, thin] (J0) -- (Ja);
\draw[black, thin] (J0) -- (Jc);
\draw[black, thin] (J0) -- (Je);
\draw[black, thin] (Ja) -- (Je);
\draw[black, thin] (Jc) -- (Je);

\draw[chain1, ->] (J0) -- (Ja);
\draw[chain1, ->] (Ja) -- (Je);

\draw[chain2, ->] (J0) -- (Jc);
\draw[chain2, ->] (Jc) -- (Je);
\draw[black, thin, ->] (J0) -- (Je);

\filldraw[black] (J0) circle (1.2pt);
\filldraw[black] (Je) circle (1.2pt);
\end{tikzpicture}
\captionsetup{labelsep=none}      
\captionof{figure}{}              
\label{fig:diamond(lingxing) Comp(S/J)}
\end{minipage}
\end{figure}

The divisibility graph $\operatorname{Div}(S)$ is shown in Figure \ref{fig:diamond(lingxing)-scc-uc-Div(S)}. The red dashed boxes are strongly connected components (SCC). The blue and orange dashed boxes are maximal unilaterally connected subgraphs $\widehat{M}_1=J_0\cup J_a\cup J_e$ and $\widehat{M}_2=J_0\cup J_c\cup J_e$. The whole graph is weakly connected (green dashed box).

\begin{figure}[htbp]
\centering
\begin{tikzpicture}[
    element/.style={circle, draw, inner sep=0pt, minimum size=0.5cm, font=\scriptsize},
    sccbox/.style={draw=red, thick, dashed, inner sep=0.3cm, rounded corners},
    ucbox/.style={draw=blue, thick, dashed, inner sep=0.5cm, rounded corners},
    >=stealth,
]

\node[element] (0) at (0,-0.5) {$0$};

\node[element] (a) at (-3.5, 0.3) {$a$};
\node[element] (b) at (-2.5, 0.4) {$b$};

\node[element] (c) at (3.5, 0.3)  {$c$};
\node[element] (d) at (2.5, 0.4)  {$d$};

\node[element] (e) at (-0.8, 2.0) {$e$};
\node[element] (f) at (0.8, 2.0)  {$f$};

\foreach \x in {a, b, c, d, e, f}
    \draw[->] (0) -- (\x);

\draw[gray, very thick, <->] (a) -- (b);
\draw[gray, very thick, <->] (c) -- (d);
\draw[gray, very thick, <->] (e) -- (f);

\foreach \x in {a, b}
    \foreach \y in {e, f}
        \draw[->] (\x) -- (\y);

\foreach \x in {c, d}
    \foreach \y in {e, f}
        \draw[->] (\x) -- (\y);

\node[sccbox, fit=(a)(b)] (scc_a) {};
\node[red, font=\footnotesize] at (scc_a.south) [below] {SCC};
\node[sccbox, fit=(c)(d)] (scc_c) {};
\node[red, font=\footnotesize] at (scc_c.south) [below] {SCC};
\node[sccbox, fit=(e)(f)] (scc_e) {};
\node[red, font=\footnotesize] at (scc_e.north) [below] {SCC};

\node[ucbox, fit=(0)(a)(b)(e)(f)] (uc1) {};
\node[blue, font=\footnotesize] at (uc1.north) [below, yshift=-2pt, xshift=-1.8cm] {\(\widehat{M}_1\) (unilateral)};
\node[ucbox, draw=orange, fit=(0)(c)(d)(e)(f)] (uc2) {};
\node[orange, font=\footnotesize]at (uc2.north) [below, yshift=-2pt, xshift=1.8cm] {\(\widehat{M}_2\) (unilateral)};

\node[draw=green, thick, dashed, inner sep=0.8cm, fit=(0)(a)(b)(c)(d)(e)(f)] (wc) {};
\node[green, font=\footnotesize] at (wc.south) [below] {Weakly connected (whole graph)};
\end{tikzpicture}
\captionsetup{labelsep=none}  
\caption{}
\label{fig:diamond(lingxing)-scc-uc-Div(S)}
\end{figure}

Let $(X_1, \le_{X_1})$ and $(X_2, \le_{X_2})$ be two posets. A map $f: X_1 \to X_2$ is called an order embedding if $f$ is injective and satisfies the condition $x\le_{X_1} y\iff f(x)\le_{X_2}f(y)$. 
An order isomorphism is a bijective order embedding, equivalently, it is a bijection such that $x\le_{X_1} y\iff f(x)\le_{X_2}f(y)$ for all $x,y\in X_1$.

The following proposition gives necessary and sufficient conditions for the weighted blow-ups of two weighted posets to be isomorphic.

\begin{proposition}[Isomorphism of weighted blow-ups]
\label{prop:weighted-blowup-isomorphism}
Let $(P,\{I_x\}_{x\in P})$ and $(Q,\{K_y\}_{y\in Q})$ be two weighted posets.  
Then 
\[
\mathsf{Bl}(P,\{I_x\}_{x\in P}) \cong \mathsf{Bl}(Q,\{K_y\}_{y\in Q})
\]
if and only if there exists a map $\varphi:P\to Q$ satisfying
\begin{enumerate}[label=(\roman*)]
  \item $\varphi$ is an order isomorphism;
  \item For any $x\in P$, $|I_x| = |K_{\varphi(x)}|$.
\end{enumerate}
\end{proposition}

\begin{proof}
\noindent\textbf{Sufficiency.}
Suppose that there exists a map $\varphi \colon P\to Q$ satisfying conditions (i) and (ii). Then $|P|=|Q|$. 
For each $x\in P$, from condition (ii), there exists a bijection
$\varphi_x \colon I_x\to K_{\varphi(x)}$.  
Hence, there exists a bijection
$\Psi\colon V\bigl(\mathsf{Bl}(P,\{I_x\}_{x\in P})\bigr)\to V\bigl(\mathsf{Bl}(Q,\{K_y\}_{y\in Q})\bigr)$, defined by $\Psi(x,\alpha)=\bigl(\varphi(x),\varphi_x(\alpha)\bigr)$.

For any distinct $(x,\alpha),(y,\beta)\in V\bigl(\mathsf{Bl}(P,\{I_x\}_{x\in P})\bigr)$,
\[
\begin{aligned}
(x,\alpha)\to(y,\beta)
&\iff x\le y\\
&\iff \varphi(x)\le \varphi(y)\\
&\iff \bigl(\varphi(x),\varphi_x(\alpha)\bigr)\to\bigl(\varphi(y),\varphi_y(\beta)\bigr)\\
&\iff \Psi(x,\alpha)\to\Psi(y,\beta).
\end{aligned}
\]
Thus, $\Psi$ is an isomorphism from $\mathsf{Bl}(P,\{I_x\}_{x\in P})$ to
$\mathsf{Bl}(Q,\{K_y\}_{y\in Q})$.

\medskip

\noindent\textbf{Necessity.}
Suppose that
$\Psi\colon\mathsf{Bl}(P,\{I_x\}_{x\in P})\to\mathsf{Bl}(Q,\{K_y\}_{y\in Q})$
is a graph isomorphism.  
In the weighted blow-up $\mathsf{Bl}(P,\{I_x\}_{x\in P})$, each set
$C_x = \{x\}\times I_x = \{(x,\alpha)\mid \alpha\in I_x\}$
is a strongly connected component.
The same description holds for
$D_y = \{y\}\times K_y = \{(y,\beta)\mid \beta\in K_y\}$
in $\mathsf{Bl}(Q,\{K_y\}_{y\in Q})$.

Since $\Psi$ is a graph isomorphism, it maps strongly connected
components to strongly connected components. Therefore, for each $x\in P$,
there exists a unique $y\in Q$ such that $\Psi(C_x)=D_y$.
This induces a map $\varphi\colon P\to Q$
defined by
$\Psi\bigl(\{x\}\times I_x\bigr)=\{\varphi(x)\}\times K_{\varphi(x)}$, i.e., $\Psi(C_x)=D_{\varphi(x)}$.
Obviously, $\varphi$ is a bijection.
For any $x\in P$, $|I_x|=|C_x|=|\Psi(C_x)|=|D_{\varphi(x)}|=|K_{\varphi(x)}|$, which gives condition (ii).

Now we show that $\varphi$ is an order isomorphism. Let $x,y\in P$ with
$x\ne y$, and choose arbitrary $\alpha\in I_x$ and $\beta\in I_y$. Since
$\Psi$ maps $\{x\}\times I_x$ to $\{\varphi(x)\}\times K_{\varphi(x)}$,
there exist $\alpha'\in K_{\varphi(x)}$ and $\beta'\in K_{\varphi(y)}$
such that
\[
\Psi(x,\alpha)=\bigl(\varphi(x),\alpha'\bigr),
\qquad
\Psi(y,\beta)=\bigl(\varphi(y),\beta'\bigr).
\]
Then
\[
\begin{aligned}
x<y
&\iff (x,\alpha)\to(y,\beta) \\
&\iff \Psi(x,\alpha)\to\Psi(y,\beta) \\
&\iff \bigl(\varphi(x),\alpha'\bigr)\to\bigl(\varphi(y),\beta'\bigr) \\
&\iff \varphi(x)<\varphi(y).
\end{aligned}
\]
Since $\varphi$ is bijective, this implies that for all $x,y\in P$,
$x\le y \iff \varphi(x)\le\varphi(y)$.
Hence $\varphi$ is an order isomorphism, which gives condition (i).
\end{proof}

By Theorem~\ref{thm: Divisibility graph blow-up theorem} and 
Proposition~\ref{prop:weighted-blowup-isomorphism}, we have the following 
corollary.

\begin{corollary}
\label{cor:isomorphism-classification}
Let $S$ and $T$ be two semigroups. 
Then $\operatorname{Div}(S)\cong\operatorname{Div}(T)$
if and only if there exists a map
$\varphi:S/\mathcal{J}\to T/\mathcal{J}$ satisfying
\begin{enumerate}[label=(\roman*)]
  \item $\varphi$ is an order isomorphism;
  \item For any $J\in S/\mathcal{J}$, $|J|=|\varphi(J)|$.
\end{enumerate}
\end{corollary}

\section{Divisibility D-saturated property of semigroups}
Having established in Section 3 that $\operatorname{Div}(S)$ is the weighted blow-up of the quotient poset $S/\mathcal{J}$, we now turn to an application of this structural description to the $D$-saturated property. 
The results of Quinn~\cite{Quinn2002} and Guo et al.~\cite{Guo2025} show that, for an arbitrary semigroup $S$, this property is determined by the quotient poset $(S/\mathcal{J},\le)$. 
It is therefore natural to ask whether the property is inherited by subsemigroups. In this section, we answer this question by introducing divisibility-preserving subsemigroups and proving an inheritance theorem for the $D$-saturated property.

\begin{lemma}
\textup{\cite{Quinn2002}}
\label{lemma: D-saturated(cycle)--equivalent}
Let $D$ be a finite graph that contains a directed cycle and $S$ be a semigroup. Then $S$ is divisibility $D$-saturated if and only if $|S/\mathcal{J}|$ is finite.
\end{lemma}

\begin{lemma}
\textup{\cite{Quinn2002}}
\label{lemma: D-saturated(acyclic)--equivalent}
Let $D$ be a finite non-null graph that is acyclic and $S$ be a semigroup. Then $S$ is divisibility $D$-saturated if and only if $S/\mathcal{J}$ has no infinite anti-chain.
\end{lemma}

Let $S$ be a semigroup and $T$ a subsemigroup of $S$. 
In the following, we write $\mathcal{J}^S$ and $\mathcal{J}^T$ for Green's $\mathcal{J}$-relations on $S$ and $T$, respectively. 
For each $u\in T$, we denote by $J_u^S$ and $J_u^T$ the $\mathcal{J}$-classes of $u$ in $S$ and $T$, and we write $S/\mathcal{J}^S$ and $T/\mathcal{J}^T$ for the corresponding quotient sets.

In general, subsemigroups do not necessarily inherit the divisibility $D$-saturated. 
For example,
let $S = B(G,\mathbb{N})$ be the Brandt semigroup over a group $G$ with an index set 
$\mathbb{N}$. Its elements are all triples $(i,g,\lambda)$ with $g\in G$, 
$i,\lambda\in\mathbb{N}$, together with a zero element $0$, and its multiplication 
is defined by
\[
(i,g,\lambda)(j,h,\mu) = 
\begin{cases}
(i,gh,\mu), & \text{if } \lambda = j,\\
0, & \text{if } \lambda \neq j.
\end{cases}
\]
As a special completely $0$-simple semigroup, $S$ has exactly two $\mathcal{J}^S$-classes, namely the non‑zero class $S\setminus\{0\}$ and the zero class $\{0\}$. Consequently $S$ is divisibility $D$-saturated for any finite graph $D$ that contains a directed cycle (Lemma~\ref{lemma: D-saturated(cycle)--equivalent}). Furthermore, $S/\mathcal{J}^S$ has no infinite anti-chain. Then $S$ is divisibility $D$-saturated for any finite graph $D$ that is acyclic (Lemma~\ref{lemma: D-saturated(acyclic)--equivalent}).

Let $T = \{(i,g,i)\mid g\in G,\; i\in\mathbb{N}\}\cup\{0\}$. It is easy to check that $T$ is an infinite subsemigroup of $S$ and $T$ has infinitely many $\mathcal{J}^T$-classes, namely $\mathcal{J}_0=\{0\},\mathcal{J}_i=\{(i,g,i)\mid g\in G\},i\in \mathbb{N}$. Consequently, $T$ is not divisibility $D$-saturated for any finite graph $D$ that contains a directed cycle. 
On the other hand, $T/\mathcal{J}^T$ has an infinite anti-chain $\{\mathcal{J}_1,\mathcal{J}_2,\cdots,\mathcal{J}_i\cdots\mid i\in \mathbb{N}\}$ since $(i,g,i)\notin T^1(j,g,j)T^1$ for any distinct $i,j\in\mathbb{N}$. Thus, $T$ is not divisibility $D$-saturated for any finite non-null graph $D$ that is acyclic.

In the following, we investigate a class of subsemigroups with the divisibility $D$-saturated property.

\begin{definition}
\label{def: divisibility-retentive}
Let $S$ be a semigroup and $T$ a subsemigroup of $S$, T is called divisibility-preserving if for all $u,v\in T$, $J_u^{S}\le J_v^{S}\;\Longleftrightarrow\; J_u^{T}\le J_v^{T}$.
\end{definition}

\begin{example}
\label{exam: divisibility-retentive}
$X=\{a,b\}$, $X^+$ is the free semigroup over $X$. Let $S=X^+$ and $T_1=\{a^n\mid n\in\mathbb{N}\}$. It is easy to check that $T_1$ is a subsemigroup of $S$. For any $x,y\in S, x\in S^1yS^1$ if and only if $y$ is a subword of $x$. For any distinct $u,v\in T_1$, there exist distinct $k,m\in \mathbb{N}$ such that $u=a^k,v=a^m$. If $J_u^{S}\le J_v^{S}$, then $u\in S^1vS^1$,i.e., $a^{k}\in S^1a^mS^1$, then $m < k$, and so $a^k\in T_1^1a^mT_1^1$,i.e., $u\in T_1^1vT_1^1$, which implies that $J_u^{T_1}\le J_v^{T_1}$. Therefore, $T_1$ is a divisibility-preserving subsemigroup of $S$. 

Let $T_2 = \{ w \in S \mid |w| \ge 2 \}$. Obviously, $T_2$ is a subsemigroup of $S$.
$T_2$ is not divisibility-preserving since $a^3\in S^1a^2S^1$, but $a^3\notin T_2^1a^2T_2^1$.
\end{example}

\begin{proposition}
\label{prop: divisibility-retentive----equivalent}
Let $S$ be a semigroup and $T$ a subsemigroup of $S$. Then the following statements are equivalent.
\begin{enumerate}
  \item[(i)] $T$ is divisibility-preserving, equivalently, for all distinct $u,v\in T$, $u\in S^1vS^1 \Longleftrightarrow u\in T^1vT^1$;
  \item[(ii)] $\operatorname{Div}(T)=\operatorname{Div}(S)|_{T}$, where $\operatorname{Div}(S)|_{T}$ denotes the subgraph of $\operatorname{Div}(S)$ induced by $T$;
  \item[(iii)] the map
        \[
        \Phi : T/{\mathcal{J}^T}\longrightarrow S/{\mathcal{J}^S},\qquad
        \Phi(J_u^{T})=J_u^{S},
        \]
        is an order embedding.
\end{enumerate}
\end{proposition}

\begin{proof}
\noindent\textbf{(i) $\Leftrightarrow$ (ii).}
Both graphs \(\operatorname{Div}(T)\) and \(\operatorname{Div}(S)|_T\) have vertex set \(T\). 
An ordered pair $(u,v)$ with $u\ne v$ is an edge of $\operatorname{Div}(T)$ if and only if $u\in T^{1}vT^{1}$, and an edge of $\operatorname{Div}(S)|_{T}$ if and only if $u\in S^{1}vS^{1}$.  
Hence, (ii) holds if and only if (i) holds.

\noindent\textbf{(i) $\Rightarrow$ (iii).}
If $J_u^{T}=J_v^{T}$, then $J_u^{S}=J_v^{S}$, i.e., $\Phi(J_u^{T})=\Phi(J_v^{T})$. Thus $\Phi$ is well defined.
Injectivity follows similarly. If $\Phi(J_u^{T})=\Phi(J_v^{T})$,
then $J_u^{S}=J_v^{S}$, which by (i) implies
$J_u^{T}=J_v^{T}$.
Furthermore, for any $u,v\in T$,
\[
J_u^{T}\le J_v^{T}
\Longleftrightarrow J_u^{S}\le J_v^{S}
\Longleftrightarrow \Phi(J_u^{T})\le \Phi(J_v^{T}),
\]
so $\Phi$ is an order embedding.

\noindent\textbf{(iii) $\Rightarrow$ (i).}
For $u,v\in T$, since $\Phi$ is an order embedding, it follows that
\[
J_u^{T}\le J_v^{T}
\Longleftrightarrow \Phi(J_u^{T})\le \Phi(J_v^{T})
\Longleftrightarrow J_u^{S}\le J_v^{S},
\]
This shows that $T$ is divisibility-preserving.
\end{proof}

\begin{lemma}
\label{lemma: divisibility-retentive subsemigroup--->D saturated(cycle)}
Let $D$ be a finite graph that contains a directed cycle and $S$ be a semigroup that is divisibility $D$-saturated. If $T$ is a divisibility-preserving subsemigroup of $S$, then $T$ is divisibility $D$-saturated.      
\end{lemma}

\begin{proof}
Suppose that $T$ is a divisibility-preserving subsemigroup of $S$. 
Since $D$ is a finite graph that contains a directed cycle and $S$ is divisibility $D$-saturated, by Lemma \ref{lemma: D-saturated(cycle)--equivalent}, it follows that $|S/\mathcal{J}^S|<\infty$. 
By Proposition \ref{prop: divisibility-retentive----equivalent}(iii), we have $|T/\mathcal{J}^T|\le |S/\mathcal{J}^S|<\infty$. 
Consequently, $T$ is divisibility $D$-saturated. 
\end{proof}

\begin{lemma}
\label{lemma: divisibility-retentive subsemigroup--->D saturated(acyclic)}
Let $D$ be a finite non-null graph that is acyclic and $S$ be a semigroup that is divisibility $D$-saturated. If $T$ is a divisibility-preserving subsemigroup of $S$, then $T$ is divisibility $D$-saturated.          
\end{lemma}

\begin{proof}
Suppose that $T$ is divisibility-preserving. 
Let $A$ be an anti-chain of $T/\mathcal{J}^T$. 
Then for any distinct $J_x^{T},J_y^{T}\in A$, $J_x^{T}$ and $J_y^{T}$ are incomparable. 
By Proposition \ref{prop: divisibility-retentive----equivalent}(iii), $\Phi(J_x^{T})$ and $\Phi(J_y^{T})$ are incomparable. 
Thus, $\Phi(A) = \{ \Phi(J_x^{T}) \mid J_x^{T} \in A\}$ is an anti-chain of $S/\mathcal{J}^S$.
If $A$ is an infinite anti-chain of $T/\mathcal{J}_T$, since $\Phi$ is injective, $\Phi(A)$ is an infinite anti-chain of $S/\mathcal{J}^S$. This contradicts the fact that $S$ is divisibility $D$-saturated.    
Hence $T/\mathcal{J}^T$ contains no infinite anti-chain.
By Lemma \ref{lemma: D-saturated(acyclic)--equivalent}, $T$ is divisibility $D$-saturated.
\end{proof}

By Lemma \ref{lemma: divisibility-retentive subsemigroup--->D saturated(cycle)} and Lemma \ref{lemma: divisibility-retentive subsemigroup--->D saturated(acyclic)}, we have that

\begin{theorem}
\label{thm: divisibility retentive subsemigroup---->D-saturated(cycle and cayclic)}
Let $D$ be a finite graph and $S$ be a semigroup that is divisibility $D$-saturated. If $T$ is a divisibility-preserving subsemigroup of $S$, then $T$ is divisibility $D$-saturated.          
\end{theorem}

Generally, the converse of Theorem~\ref{thm: divisibility retentive subsemigroup---->D-saturated(cycle and cayclic)} is not true. 
Let $D$ be a finite graph and $S = B(G,\mathbb{N})$. So $S$ is divisibility $D$-saturated. 
Indeed, the subsemigroup $T_3=\{(1,g,i)\mid g\in G,i=1,2\}\cup\{0\}$ is also divisibility $D$-saturated since $T_3$ has only three $\mathcal{J}_{T_3}$-classes. But it is not divisibility-preserving since $(1,g,1)\in S^1(1,h,2)S^1$ but $(1,g,1)\notin T_3^1(1,h,2)T_3^1$.

\vspace{0.3cm}   
\noindent\textbf{Acknowledgements } The authors are particularly grateful to Professor Mikhail V. Volkov for his valuable comments and suggestions, which significantly improved the paper. Xuan Yang, Yong Shao, and Lingli Zeng are supported by Shaanxi Provincial Natural Science Basic Research Program Project (No. 2025JC-YBMS-086).

\end{document}